\documentclass[12pt,twoside]{amsart}
\usepackage{amssymb}
\usepackage{amscd}

\title[Algebraic fiber spaces]
{Algebraic fiber spaces whose general fibers are of 
maximal Albanese dimension}
\author{Osamu Fujino}
\date{2002/4/22}
\address{Research Institute for Mathematical Sciences\\ 
Kyoto University, Kyoto 606-8502 Japan}
\email{fujino@kurims.kyoto-u.ac.jp}
\newcommand{\bQ}[0]{{\mathbb Q}}

\newcommand{\Supp}[0]{{\operatorname{Supp}}}
\newcommand{\xIm}[0]{{\operatorname{Im}}}

\newcommand{\Alb}[0]{{\operatorname{Alb}}}
\newcommand{\alb}[0]{{\operatorname{alb}}}
\newcommand{\xExc}[0]{{\operatorname{Exc}}}
\newcommand{\xVar}[0]{{\operatorname{Var}}}
\newcommand{\Spec}[0]{{\operatorname{Spec}}}
\newcommand{\codim}[0]{{\operatorname{codim}}}

\newtheorem{thm}{Theorem}[section]
\newtheorem{lem}[thm]{Lemma}
\newtheorem{cor}[thm]{Corollary}
\newtheorem{prop}[thm]{Proposition}
\newtheorem{conj}[thm]{Conjecture}

\theoremstyle{definition}
\newtheorem{defn}[thm]{Definition}
\newtheorem{ex}[thm]{Example}
\newtheorem{rem}[thm]{Remark}
\newtheorem*{ack}{Acknowledgements}       

\newtheorem*{notation}{Notation}         
       
\newtheorem{say}[thm]{}

\theoremstyle{remark}

\begin{document}
\bibliographystyle{amsalpha+}

\abstract 
The main purpose of this paper is 
to prove the Iitaka conjecture $C_{n,m}$ on 
the assumption that 
the sufficiently general fibers have maximal Albanese dimension. 
\endabstract
\maketitle

\setcounter{section}{-1}
\section{Introduction}\label{se1}

The main purpose of this paper is 
to prove the Iitaka conjecture $C_{n,m}$ on 
the assumption that 
the sufficiently general fibers have maximal Albanese dimension. 
For the definition of varieties of maximal Albanese dimension, 
see Definition \ref{mad} below. 
The proof of the main theorem was essentially given 
in \cite[Section 7]{fu}, which is a variant of 
\cite[Proof of Theorem 16]{Ka}. 

\vspace{5mm}
The following is the main theorem of this paper. 

\begin{thm}\label{int}
Let $f:X\to Y$ be a proper surjective morphism 
between non-singular projective varieties with connected fibers. 
Assume that sufficiently general fibers are of 
maximal Albanese dimension. 
Then $\kappa (X)\geq \kappa (Y)+\kappa 
(X_{\eta})$, where $X_\eta$ is the generic 
fiber of $f$. 
\end{thm}

Iitaka's conjecture $C_{n,m}$ was proved on the assumption 
that the general fibers are of general type by 
Koll\'ar. For the details, see \cite{ko} or Theorem \ref{ge}. 
We use Koll\'ar's theorem to prove Theorem \ref{int}. 
By the definition of general types, 
pluricanonical maps are birational. 
So, pluricanonical maps don't lose birational properties 
of varieties of general type. 
On the other hand, 
varieties of maximal Albanese dimension have Albanese maps 
that are generically finite on their images. 
Therefore, Albanese maps lose properties of varieties 
of maximal Albanese dimension little. 
So, it is not surprising that we can prove the above theorem. 

We summarize the contents of this paper: 
Section \ref{2} contains preliminaries. 
We recall Iitaka's conjecture $C_{n,m}$ and 
some known results. 
In Section \ref{mad3}, we define varieties of {\em{maximal 
Albanese dimension}} and collect some basic properties of them. 
Section \ref{3} deals with a canonical bundle formula, 
which was obtained in \cite{fm} and \cite{fu},  
{\em{semistable parts}} and so on. We investigate the 
relationship between semistable parts and 
{\em{variations}}. 
In Section \ref{4}, we will prove the main theorem. 
As stated above, the essential part of the proof was 
contained in \cite{Ka}.   
In Section \ref{5}, we treat some results about Abelian varieties 
which we need in the 
proof of the main theorem.   

\begin{ack} 
I was inspired by the preprint \cite{ch}. 
Some parts of this paper was done during the visit to 
Newton Institute in University of Cambridge. 
I am grateful to the institute for providing 
an excellent working environment. 
I was partially supported by 
Inoue Foundation for Science. 
I am grateful to Professor Shigeru Mukai and Doctor 
Hokuto Uehara for 
giving me some comments. 
I like to thank Professor Noboru Nakayama for 
telling me Example \ref{tuika} (2). 
\end{ack}

We fix the notation used in this paper. 

\begin{notation}\label{no} 
We will work over the complex number field $\mathbb C$ 
throughout this paper. 
\begin{itemize}
\item[(i)] 
A {\em{sufficiently general point}} $z$ (resp.~{\em{subvariety}} 
$\Gamma$) of the variety $Z$ 
means that $z$ (resp.~$\Gamma$) is not contained in the countable union of 
certain proper Zariski closed subsets. 

Let $f:X\to Y$ be a morphism between varieties. 
A {\em{sufficiently general fiber}} $X_y=f^{-1}(y)$ of $f$ means 
that $y$ is a sufficiently general point in $Y$. 

\item[(ii)] 
An {\em{algebraic fiber space}} $f:X\to Y$ is a proper surjective 
morphism between non-singular projective varieties 
$X$ and $Y$ with connected fibers. 

\item[(iii)]
Let $X$ be a smooth projective variety. If 
the Kodaira dimension $\kappa(X)>0$,
then we have the Iitaka fibration $f: X'\to Y$, 
where $X'$ and $Y$ are non-singular 
projective varieties, $X'$ is birationally equivalent to $X$, and
$Y$ is of dimension
$\kappa(X)$, such that the sufficiently 
general fiber of $f$ is smooth,
irreducible with $\kappa =0$.  
The Iitaka fibration is determined only up to
birational equivalence.
Since we are interested in questions of a birational nature,
we usually assume that
$X=X'$ and that $Y$ is smooth.
Note that we often modify $f:X\to Y$ birationally without mentioning it. 
For the basic properties of the Kodaira dimension and 
the Iitaka fibration, 
see \cite[Chapter III]{ueno} or \cite[Sections 1,2]{M}. 

\item[(iv)] 
Let $B_+, B_-$ be the effective $\mathbb Q$-divisors on a variety $X$
without common irreducible components such that 
$B_+-B_-=B$. They are called the {\em {positive}}
and the {\em{negative}} parts of $B$.

Let $f : X \to Y$ be a surjective morphism.
Let $B^h, B^v$ be the $\mathbb Q$-divisors on $X$
with $B^h+B^v=B$
such that an irreducible component of $\Supp B$
is contained in $\Supp B^h$ if and only if it is mapped onto $Y$.
They are called the {\em horizontal}
and the {\em vertical} parts of $B$ over $Y$.
A divisor $B$ is said to be {\em horizontal} (resp.~{\em vertical}\ )
over $Y$ if $B=B^h$ (resp.~$B=B^v$). 
The phrase ``over $Y$" might be suppressed
if there is no danger of confusion.

\item[(v)] 
Let $\varphi:V\to W$ be a generically finite morphism 
between varieties. 
By the {\em{exceptional locus}} of $\varphi$, 
we mean the subset $\{v\in V | \dim \varphi^{-1}\varphi(v)\geq 1\}$ 
of $V$, and denote it by $\xExc(\varphi)$. 
\end{itemize}
\end{notation}

\section{Preliminaries}\label{2}

We recall the Iitaka conjecture. 
The following is a famous conjecture by Iitaka \cite [p.26 Conjecture C]
{i}. 
For the details, see \cite [Sections 6, 7]{M}. 

\begin{conj}[{Conjecture $C_{n,m}$}]\label{conj}
Let $f:X\to Y $ be an algebraic fiber space with 
$\dim X=n$ and $\dim Y= m$. 
Then we have 
$$
\kappa (X)\geq \kappa (Y)+\kappa (X_{\eta}),  
$$
where $X_\eta$ is the generic fiber of $f$. 
\end{conj}

We recall some known results about the above conjecture, 
which will be used in the proof of the main theorem. 
The following is a part of \cite[p.363 Theorem]{ko}. 
We note that a simplified proof was obtained by Viehweg (see \cite
[Theorem 1.20]{v2}). 

\begin{thm}\label{ge} 
Let $f:X\to Y$ be an algebraic fiber space such that 
the generic fiber of $f$ is of general type. 
Then we have; 
$$
\kappa (X)\geq \kappa (Y)+\kappa (X_{\eta}). 
$$
\end{thm}

Let us recall the notion of {\em{variation}} (cf.~\cite[p.329]{v}). 

\begin{defn}\label{var} 
Let $f:X\to Y$ be an algebraic fiber space. 
The {\em{variation}} of $f$, which is denoted by  
$\xVar(f)$, is defined to be the minimal number $k$, 
such that there exists a subfield $L$ of 
$\overline{\mathbb C(Y)}$ of transcendental degree 
$k$ over $\mathbb C$ and a variety $F$ over $L$ 
with $F\times _{\Spec (L)}\Spec ({\overline {\mathbb C(Y)}})\sim 
X \times_Y \Spec (\overline{{\mathbb C}(Y)})$, where 
$\sim$ means ``birational''. 
\end{defn}

The next theorem is in \cite[Corollary 14]{Ka}. It is also a special 
case of \cite[Corollary 1.2 (ii)]{ka3}. 
See also \cite[Corollary 7.3]{fu}. 

\begin{thm}\label{takusan}
Let $f:X\to Y$ be an algebraic fiber space such 
that the geometric generic fiber is birationally equivalent 
to an Abelian variety. 
If $\kappa (Y)\geq 0$, then we have 
$$
\kappa (X)\geq \max\{\kappa (Y), \xVar (f)\}.  
$$
\end{thm}

\section{Varieties of maximal Albanese dimension}\label{mad3}

Let us recall the definition of the varieties of 
maximal Albanese dimension. I learned it from \cite {ch} 
and \cite{HP}. 

\begin{defn}[Varieties of maximal Albanese dimension]
\label{mad}
Let $X$ be a smooth projective variety. 
Let $\Alb (X)$ be the Albanese variety 
of $X$ and $\alb_X :X\to \Alb (X)$ the 
corresponding Albanese map. 
We say that 
$X$ has {\em{maximal Albanese dimension}}, or 
is {\em{of maximal Albanese dimension}}, if 
$\dim (\alb _X(X))=\dim X$. 
\end{defn}

\begin{rem}\label{doudemo}
A smooth projective variety $X$ has maximal Albanese dimension 
if and only if the cotangent bundle of $X$ is 
generically generated by its global sections, 
that is, 
$$
H^{0}(X, \Omega_X^{1})\otimes_{\mathbb C}\mathcal O_X\to 
\Omega _X^{1} 
$$ 
is surjective at the generic point of $X$. 
It can be checked without any difficulty. 
\end{rem} 

For the basic 
properties of Albanese mappings, see \cite[Chapter IV \S 9]
{ueno}. 

\begin{prop}\label{kote}
The following properties are easy to check by the definition. 
\begin{itemize}
\item[(1)] The notion of maximal Albanese dimension is 
birationally invariant. 
\item[(2)] Let $X$ be an Abelian variety. 
Then $X\simeq \Alb(X)$. 
Of course, it has maximal Albanese dimension. 
\item[(3)] Let $X$ be a variety of maximal Albanese dimension. 
Let $Y$ be a smooth projective variety and 
$\varphi:Y\to X$ a morphism such that 
$\dim Y=\dim \varphi (Y)$. 
If $\varphi(Y)\not\subset\xExc(\alb_X)$, 
then $Y$ has maximal Albanese dimension. 
\proof[Proof of (3)]
By the universality of Albanese mappings, 
we have the following commutative diagram: 
$$
\begin{CD}
Y &@>{\alb _Y}>> & \Alb (Y)\\
@V{\varphi}VV& & @V{\varphi_*}VV \\
X &@>{\alb _X}>> & \Alb (X).
\end{CD}
$$
So, we have that $\varphi_*(\alb_Y(Y))=\alb_X(\varphi(Y))$. 
By the assumption, $\dim (\alb_X(\varphi(Y)))=\dim Y$ holds. 
Therefore, we have that $\dim(\alb_Y(Y))=\dim Y$. 
This means that $Y$ has maximal Albanese dimension. 
\endproof
\end{itemize}
\begin{itemize}
\item[(4)] Let $X$ and $Y$ be varieties of 
maximal Albanese dimension. 
Then, so is $X\times Y$. 
It is obvious since $\Alb(X\times Y)\simeq\Alb(X)\times
\Alb(Y)$. 
\item[(5)] Let $X$ be a variety 
of maximal Albanese dimension. 
Then the Kodaira dimension $\kappa (X)\geq 0$. 
If $\kappa (X)=0$, then $X$ is birationally equivalent to 
$\Alb(X)$ by \cite[Theorem 1]{ka0}. 
\end{itemize}
\end{prop}

The following proposition is \cite[Proposition 2.1]{HP}. 
See also \cite[Theorem 13]{ka0}. 

\begin{prop}\label{albanese}
Let $X$ be a smooth projective variety of maximal
Albanese dimension, and let
$f\colon X\to Y$ be the Iitaka fibration. 
We can assume that $Y$ is smooth by 
using Hironaka's desingularizaion theorem {\em{(}}
see Notation (iii) 
and Proposition \ref{kote} (1){\em{)}}. 
We have the following commutative diagram by the universal 
property 
of Albanese varieties: 
$$
\begin{CD}
X &@>{\alb _X}>> & \Alb (X)\\
@V{f}VV& & @V{f_*}VV \\
Y &@>{\alb _Y}>> & \Alb (Y).\\
\end{CD}
$$
Then we have: 
\begin{itemize}
\item[(a)] $Y$ has maximal Albanese dimension;
\item[(b)] $f_*$ is surjective and $\ker f_*$ is connected, of 
dimension $\dim (X)-\kappa(X)$;
\item[(c)] there exists an Abelian variety
$P$ isogenous to $\ker f_*$ such that the 
sufficiently general fiber of $f$ is birationally equivalent to
$P$.
\end{itemize}
\end{prop}
\proof[Sketch of the proof] 
By Proposition \ref{kote} (3) and (5), we can check that 
the sufficiently general fibers of $f$ are birationally equivalent to 
an Abelian variety. By easy dimension count, we can prove (a) and 
(b) without difficulty. 
For the details, see \cite[Proposition 2.1]{HP}. 
\endproof

The next example helps the readers to understand 
Proposition \ref{albanese}. 

\begin{ex}\label{tuika}
\begin{enumerate}
\item[(1)] Let $C$ be a smooth projective curve 
with genus $g(C)\geq 1$. Then $\alb_C:C\to \Alb(C)$ 
is an embedding by Abel's theorem. 
Therefore, $C$ has maximal Albanese dimension. 
\item[(2)] Let $C$ be a smooth projective curve with 
$g(C)\geq 2$. 
We assume that there exists an involution $\iota_1$ on 
$C$ such that $Y:=C/\iota_1$ is an elliptic curve. 
Let $E$ be an elliptic curve and 
$a\in E$ a non-zero $2$-torsion point. 
We define $\iota_2:=T_a:E\to E$, 
where $\iota_2(b)=b+a$ for every $b\in E$. 
Then $E':=E/\iota_2$ is an elliptic curve. 
We put $S:=C\times E$. 
It is obvious that $\kappa (S)=1$ and 
$S$ has maximal Albanese dimension 
by Proposition \ref{kote} (2), (4) and 
Example \ref{tuika} (1). 
Let $\iota:=\iota_1\times \iota_2:S\to S$. 
Then $\iota$ is an involution. 
We put $X:=S/\iota$. 
Since the action is free, 
$X$ is smooth and $\kappa (X)=1$. 
We can check that the projection $f:X\to Y$ 
is the Iitaka fibration. 
We note that $\kappa (X)=1$ and general fibers of $f$ 
are isomorphic to $E$. 
By the following commutative diagram, 
$$
\begin{CD}
E   & @>>> & E' & @>{id}>> E'=\Alb(E')\\
@AAA& & @AAA & @AAA\\
S   &@>>> & X &@>{\alb_X}>> \Alb(X)\\
@VVV& & @VV{f}V & @VVV\\
C   & @>>> & Y & @>{id}>> Y=\Alb(Y).
\end{CD}
$$ 
we can check that $S\to X\to \Alb(X)$ is a finite morphism. 
In particular, $\alb_X:X\to \Alb(X)$ is finite. 
Therefore, $X$ has maximal Albanese dimension and 
$\kappa (X)=1$ such that $f:X\to Y$ is the Iitaka fibration, 
where $Y$ is an elliptic curve. 
\end{enumerate}
\end{ex}

\section{Semistable parts and Variations}\label{3}

We review the basic definitions and properties 
of the semistable part $L_{X/Y}^{ss}$ without proof. 
For the details, we recommend the reader to see 
\cite [Sections 2, 4]{fm} and \cite[Sections 3, 4]{fu}. 

\begin{say}\label{(2.1)} Let $f : X \to Y$ be
an algebraic fiber space 
such that the Kodaira dimension of the generic fiber 
of $f$ is zero, that is, $\kappa(X_\eta) = 0$.
We fix the smallest $b \in \mathbb N$
such that the $b$-th plurigenus
$P_b(X_{\eta})$ is non-zero.
\end{say}

\begin{prop}[{\cite[Proposition 2.2]{fm}}]\label{(2.2)} 
There exists one and only one
$\mathbb Q$-divisor $D$ modulo linear equivalence on $Y$ with a
graded $\mathcal O_Y$-algebra isomorphism
$$
\bigoplus_{i \ge 0} \mathcal O_Y(\lfloor iD \rfloor)
\simeq \bigoplus_{i \ge 0} (f_*\mathcal O_X(ibK_{X/Y}))^{**},
$$
where $M^{**}$ denotes the double dual of $M$.

Furthermore, the above isomorphism induces the equality
$$
bK_X = f^*(b K_Y+D)+B,
$$
where $B$ is a $\mathbb Q$-divisor on $X$ such that
$f_*\mathcal O_X(\lfloor iB_+\rfloor)\simeq\mathcal O_Y\ (\forall i>0)$
and $\codim_Y f(\Supp B_-) \ge 2$.
We note that
for an arbitrary open set $U$ of $Y$, $D|_U$ and $B|_{f^{-1}(U)}$
depend only on $f|_{f^{-1}(U)}$.

If furthermore $b=1$ and fibers of $f$ over codimension one 
points of $Y$ are all reduced, then the divisor $D$ is a Weil 
divisor. 
\end{prop}

\begin{defn}\label{(2.4)}
Under the notation of \ref{(2.2)}, we denote $D$ by $L_{X/Y}$.
It is obvious that $L_{X/Y}$ depends only on the birational
equivalence class
of $X$ over $Y$. 
\end{defn}

The following definition is a special 
case of \cite [Definition 4.2]{fm} 
(see also \cite[Proposition 4.7]{fm}).

\begin{defn}\label{(x2.3)}
We set $s_P:=b(1-t_P)$, where
$t_P$ is the log-canonical threshold of ${f}^{*}P$ with
respect to $(X, -(1/b)B)$ over the
generic point $\eta_P$ of $P$:
$$
t_P:=\max\{t\in \mathbb R\ |\
(X, -(1/b)B+tf^{*}P)
\text{ is log-canonical over } \eta_P\}.
$$
Note that $t_P\in \mathbb Q$ and that $s_P\ne 0$ only for 
a finite number of codimension one points $P$
because there exists a nonempty Zariski open set $U\subset Y$ 
such that $s_P=0$ for every prime divisor $P$ with 
$P\cap U\ne \emptyset$. 
We note that $s_P$ depends only on
$f|_{f^{-1}(U)}$ where
$U$ is an open set containing $P$.

We set $L_{X/Y}^{ss}:=L_{X/Y}-\sum_P 
s_P P$
and call it the {\it semistable part} of $K_{X/Y}$.

We note that $D$, $L_{X/Y}$, $s_P$,
$t_P$ and $L_{X/Y}^{ss}$
are birational invariants of $X$ over $Y$. 

Putting the above symbols together, we have 
{\it the canonical bundle formula}
for $X$ over $Y$:
$$
bK_X = f^*(bK_Y+L_{X/Y}^{ss})+
\sum_P s_P f^*P
+ B,
$$
where $B$ is a $\bQ$-divisor on $X$ such that 
$f_*\mathcal O_X(\lfloor iB_+\rfloor)\simeq\mathcal O_Y\ 
(\forall i>0)$ and
$\codim_Y f(\Supp B_-) \ge 2$.
\end{defn}

\begin{defn}[Canonical cover of the generic fiber]\label{(2.8.1)}
Under the notation of \ref{(2.1)},
consider the following construction.
Since $\dim |bK_{X_\eta}| = 0$, there exists a Weil divisor
$W$ on $X$ such that

\begin{enumerate}
\item[(i)] $W^h$ is effective 
and $f_*\mathcal O_{X}(iW^h) \simeq\mathcal O_Y$ for all $i > 0$, and
\item[(ii)] $bK_X - W$ is a principal divisor $(\psi)$ for
some non-zero rational function $\psi$ on $X$.
\end{enumerate}

Let $s:Z \to X$ be the
normalization of $X$ in $\mathbb C(X)(\psi^{1/b})$. 
We call $Z\to X\to Y$ a {\em{canonical cover}} of $X\to Y$. 
We often call $Z'\to X$ a canonical 
cover after replacing $Z$ with its resolution $Z'$. 
\end{defn}

By using \ref{(2.8.1)} and replacing $f:X\to Y$ birationally, 
we always make the situation as in \ref{lem3.1} (see 
\cite[(5.15.2)]{M}). 

\begin{say}\label{lem3.1}
Let $f:X\to Y$ and $h:W\to Y$ be algebraic fiber spaces such 
that 
\begin{itemize}
\item[(i)] the algebraic fiber space $f$ is as in \ref{(2.1)},
\item[(ii)] $h$ factors as 
$$
\begin{CD}
h:W@>{g}>>X@>{f}>>Y, 
\end{CD}
$$
where $g$ is generically finite,  
\item[(iii)] there is a simple normal crossing divisor 
$\Sigma$ on $Y$ such that $f$ and $h$ 
are smooth over $Y_{0}:=Y\setminus \Sigma$, 
\item[(iv)] the Kodaira dimension of 
the generic fiber $W_{\eta}$ is zero and the geometric genus 
$p_{g}(W_{\eta})=1$, where $\eta$ is the generic point of 
$Y$. 
\end{itemize}
\end{say}
\begin{say}\label{3.3}
By the definition of $\xVar (f)$, there are an 
algebraic fiber space $f':X'\to Y'$ with $\overline {\mathbb C(Y')}=L$, 
a generically finite and generically surjective morphism 
$\pi:\overline Y\to Y$ and a generically surjective morphism 
$\rho:\overline Y\to Y'$ such that 
the induced algebraic fiber space 
$\overline f:\overline X\to \overline Y$ 
from $f$ by $\pi$ is birationally equivalent to 
that from $f'$ by $\rho$ as in the 
following commutative diagram: 
$$
\begin{CD}
X@<<<{\overline X}@>>>{X'}\\
@V{f}VV @VV{\overline f}V @VV{f'}V\\
Y@<<{\pi}<{\overline Y}@>>{\rho}>{Y'}. 
\end{CD}
$$
Furthermore, we can assume that $Y'$ and $\overline Y$ are smooth 
projective varieties. 
We also assume that there are simple normal crossing divisors 
$\Sigma$ on $Y$ and $\Sigma'$ on $Y'$ as in \ref{lem3.1} 
such that both $\pi^{-1}(\Sigma)$ and $\rho^{-1}(\Sigma')$ are 
simple normal crossing divisors on $\overline Y$. 
Then by \cite[Proposition 4.2]{fu}, we obtain  
that 
$$
\pi^{*}L_{X/Y}^{ss}=L_{\overline X/ \overline Y}^{ss}
=\rho^{*}L_{X'/Y'}^{ss}. 
$$ 
Therefore, we have: 
\begin{thm}\label{sin} 
Let $f:X\to Y$ be an algebraic fiber space as in \ref{lem3.1}. 
Then 
$$
\kappa (Y,L_{X/Y}^{ss})
=\kappa (Y', L_{X'/Y'}^{ss})\leq \dim Y'=\xVar (f). 
$$
\end{thm}
\end{say}

\begin{rem}
By \cite[Theorem 1.1]{ka3}, $\kappa (Y,L_{X/Y}^{ss})=\xVar (f)$ 
if there exists a good minimal algebraic variety 
$X_{{\overline \eta}, {\text{min}}}$ defined over 
$\overline{\mathbb C(Y)}$ that is birationally 
equivalent to the geometric generic 
fiber $X_{{\overline \eta}}$ over 
$\overline{\mathbb C(Y)}$. For the details, see \cite{ka3}. 
\end{rem}

\begin{say}\label{3.4}
In the same situation as in \ref{lem3.1}, 
we further assume that $\xVar(f)=0$. 
In this case, $Y'$ is a point and 
$\overline X$ is birationally equivalent to $\overline Y\times F$ 
for a non-singular projective variety $F$. 
Let $\widetilde F$ be a {\em{canonical cover}} of $F$ 
(see Definition \ref{(2.8.1)}). 
We note that we applied Definition \ref{(2.8.1)} for 
$F\to \Spec\, \mathbb C$. 
Then $$bL_{\overline Y\times \widetilde F/\overline Y}^{ss} 
=L_{\overline Y\times F/\overline Y}^{ss}$$ by 
\cite[(5.15.8)]{M} or \cite[Lemma 4.1]{fu}, where 
$b$ is the smallest positive integer 
such that the $b$-th plurigenus $P_b(F)\ne 0$. 
On the other hand, we can check easily that  
$$
L_{\overline Y\times \widetilde F/\overline Y}^{ss}\sim 0. 
$$
Thus, $mL_{X/Y}^{ss}\sim 0$ for some positive integer $m$. 
Therefore, we summarize; 
\begin{thm}\label{dou}
Let $f:X\to Y$ be an algebraic fiber space as in 
\ref{lem3.1}. Assume that $\xVar(f)=0$. Then we obtain;  
$$mL_{X/Y}^{ss}\sim 0$$ for some positive integer $m$.  
\end{thm}
\end{say}

\section{Main Theorem}\label{4}

The following is the main theorem of 
this paper. This says that 
the Iitaka conjecture $C_{n,m}$ is true on 
the assumption that the sufficiently general fibers have 
maximal Albanese dimension. 

\begin{thm}\label{main}
Let $f:X\to Y$ be an algebraic fiber space. 
Assume that sufficiently general fibers are of 
maximal Albanese dimension. 
Then $\kappa (X)\geq \kappa (Y)+\kappa 
(X_{\eta})$, where $X_\eta$ is the generic 
fiber of $f$. 
\end{thm}

\proof[Proof of the theorem]
If $\kappa (Y)=-\infty$, then the inequality 
is obviously true. 
So, we can assume that $\kappa (Y)\geq 0$. 

If $\kappa (X_\eta)=0$, then the geometric generic fiber 
is birationally equivalent to an 
Abelian variety by Proposition \ref{kote} (5). 
Thus $\kappa (X)\geq \kappa (Y)=\kappa (Y)
+\kappa (X_\eta)$ by Theorem \ref{takusan}. 
Therefore, we can assume that $\kappa (X_\eta) >0$ 
from now on. 

The following lemmas \cite[Lecture 4]{i2} are useful. 
We write it for the reader's convenience 
(see also \cite[Proposition 6]
{Ka}). 

\begin{lem}[Induction Lemma]\label{induction}
Under the same notation as in Theorem \ref{main}, 
it is sufficient to prove that 
$\kappa (X)>0$ on the assumption 
that 
$\kappa (Y)\geq 0$ and $\kappa (X_\eta)>0$. 
\end{lem}

\proof[Proof of the lemma] 
We use the induction on the dimension of $X$. 
If $\dim X=1$, then there is nothing to prove. 
 
Let $\varphi:X\to Z$ be the Iitaka fibration 
associated to $X$. 
Since $\kappa (X)>0$, we have $\dim Z=\kappa (X)>0$. 
For a sufficiently general point $z\in Z$, 
the fiber $X_z=\varphi^{-1}(z)$ has 
Kodaira dimension zero. 
We define $f':=f|_{X_z}:X_z\to B=f(X_z)$. 
We note that the sufficiently general 
fiber of $f'$ is of maximal Albanese dimension. 
By induction hypothesis, 
$$
0=\kappa (X_z)\geq \kappa (X_{z,y})+\kappa (B),  
$$ 
where $y$ is a sufficiently general point of $B$. 
By Lemma \ref{iitaka} below, 
we have $\Gamma$, $\widetilde W$, $F$ and 
$G$. 
Since $\Gamma$ is sufficiently general, we can assume that 
$z\in \Gamma$. 
Furthermore, since $\dim \widetilde W=\dim Y$, 
$\kappa (Y)\leq \kappa (\widetilde W)$ follows. 
And by the easy addition, 
we get 
\begin{eqnarray*}
\kappa(\widetilde W)&\leq& \kappa (F^{-1}(z))+\dim \Gamma\\
&=&\kappa (B)+\dim Y -\dim B. 
\end{eqnarray*}
By hypothesis, 
$\kappa (Y)\geq 0$ and hence
$$
0\leq \kappa (Y)\leq \kappa (\widetilde W)
\leq \kappa (B)+\dim Y-\dim B. 
$$
This implies $\kappa (B)\geq 0$. 

On the other hand, $X_{z,y}=f^{-1}(y)\cap\varphi^{-1}(z)$ 
can be considered as a 
sufficiently general fiber of 
$\varphi|_{X_y}:X_y\to \varphi(X_y)$, 
where $y$ is also a sufficiently general 
point of $Y$. Thus, $\kappa (X_{z,y})\geq 0$. 
More precisely, $X_{z,y}$ is of maximal Albanese dimension. 
Therefore, 
we get $\kappa (X_{z,y})=\kappa (B)=0$. 
By the easy addition, 
$$
\kappa (X_y)\leq \kappa (X_{z,y})+\dim(\varphi(X_y)). 
$$
So, $\kappa (X_y)\leq \dim (\varphi(X_y))$. 
Clearly, we have 
\begin{eqnarray*}
\dim (\varphi(X_y))&=&\dim X_y-\dim X_{z,y}\\
&=&\dim X-\dim Y-(\dim X_z-\dim B)\\
&=&\dim Z+\dim B-\dim Y\\
&=&\kappa (X)+\dim B-\dim Y.
\end{eqnarray*}
Hence, 
\begin{eqnarray*}
\kappa (X)&\geq &\kappa (X_y)+\dim Y-\dim B\\
&\geq &\kappa (Y)+\kappa (X_y). 
\end{eqnarray*}
We note that $\kappa (X_\eta)=\kappa (X_y)$. 
We finish the proof of the lemma. 
\endproof

The next lemma was already used in the proof of 
Lemma \ref{induction}. See \cite[p.46]{i2}. 

\begin{lem}[Kawamata]\label{iitaka}
Let $f:X\to Y$ and $\varphi:X\to Z$ be 
proper surjective morphisms with connected fibers, 
where $X$, $Y$, and $Z$ are normal projective varieties. 
Then there exists a sufficiently general 
subvariety 
$\Gamma$ of $Z$, a variety $\widetilde W$ 
and morphisms 
$F:\widetilde W\to \Gamma$, 
$G:\widetilde W\to Y$ such that 
$F:\widetilde W\to \Gamma$ is a proper surjective morphisms with 
$F^{-1}(z)=f(\varphi^{-1}(z))$, 
and $G:\widetilde W\to Y$ is generically finite. 
\end{lem}
\proof
Let $\Phi:=(f,\varphi):X\to Y\times Z$ and $V$ be the closure of 
$\xIm \Phi$. 
Restricting the projection morphisms, we have $p:V\to Y$ and 
$q:V\to Z$. 
For $z\in Z$, $q^{-1}(z)=(f(\varphi^{-1}(z)), z)\simeq 
f(\varphi^{-1}(z))$, 
and for $y\in Y$, $p^{-1}(y)=(y,\varphi(f^{-1}(y)))\simeq 
\varphi(f^{-1}(y))$. 
Hence $p$ is surjective and let $r=\dim V-\dim Y$. 
If $r=0$, then $\widetilde W:=V$ has the required property. 
If $r>0$, take a sufficiently general hyperplane section $Z_1$ of $Z$. 
$V_1=q^{-1}(Z_1)$ is isomorphic to $V\times _Z Z_1$ and 
also $V_1=V\cap (Y\times Z_1)$. 
Then $p_1:=p|_{V_1}:V_1\to Y$ satisfies that 
$V_1$ is a variety and $p_1^{-1}(y)=(y, \varphi(f^{-1}(y))\cap Z_1)$. 
Since $Z_1$ is sufficiently general, 
for a general point $y\in Y$, it follows that 
$\dim (\varphi(f^{-1}(y))\cap Z_1)=r-1$. 
Repeating this $r$ times, we have $\Gamma :=Z_r$ and 
$\widetilde W=V_r$ have the required property. 
\endproof

\proof[Proof of the theorem continued]
We use the induction with respect to $\dim X$ to prove the 
main theorem. 
If the generic fiber $X_{\eta}$ is of general type, 
then $\kappa (X)\geq \kappa (Y)+\dim X_{\eta}>0$ by 
Theorem \ref{ge}. So we can assume that 
$\kappa (X_{\eta})<\dim X_{\eta}$. 
Let $X\to Z\to Y$ be the relative Iitaka fibration. 
Then the geometric generic fiber of 
$g:X\to Z$ is birationally equivalent to 
an Abelian variety. 
So, $\kappa (X)\geq \max\{\kappa (Z), \xVar(g)\}$ 
by Theorem \ref{takusan}. 
We note that the sufficiently general fiber of 
$h:Z\to Y$ has maximal Albanese dimension by Proposition \ref{albanese} (a).
Therefore, $\kappa (Z)\geq 0$ by the induction 
and we can apply Theorem \ref{takusan} to 
$g:X\to Z$. 
If the Kodaira dimension of the sufficiently general fiber 
of $h$ is positive, 
then $\kappa (Z)>0$ by the induction. 
Thus we have $\kappa (X)\geq \kappa (Z)>0$. 
So, we can assume that the geometric generic fiber of 
$h$ is of Kodaira dimension zero. 
Therefore, $\kappa (Z)\geq \xVar(h)$ since 
the geometric generic fiber is birationally equivalent to 
an Abelian variety (see Theorem \ref{takusan}). 

Thus, we can assume that $\xVar(h)=\xVar(g)=0$ and 
the geometric generic fiber of $h$ is birationally equivalent to 
an Abelian variety. 

We shall prove that $\kappa (X,K_{X/Y})>0$ for 
a suitable birational
model of $f:X\to Y$. 
Using \cite[Lemma 7.8]{fu} and \cite[Theorems 8, 9]{Ka}, 
we reduce it to the case where $Z$ is birationally equivalent to 
a product $Y\times A$ for an Abelian variety $A$. 
Thus we come to the following situation: 
$$
\begin{CD}
f:X@>{g}>>Z@>{\nu}>>{Y\times A}@>{h_{1}}>>Y, 
\end{CD}
$$ 
where 
\begin{itemize}
\item[(a)] $A$ is an Abelian variety, 
\item[(b)] $f$ is the given fiber space, 
$h_1$ is the projection, and $\nu$ is a proper birational morphism, 
\item[(c)] there is a simple normal crossing divisor $D$ on 
$Z$ such that $g$ is smooth over $Z\setminus D$, 
and $f$ factors as 
$$
\begin{CD}
X@>{\mu}>>{\widetilde X}@>>>Y, 
\end{CD}
$$
where $\mu$ is birational and $\widetilde X$ is a non-singular 
projective variety 
such that $B_-$ is an effective $\mu$-exceptional divisor 
by \cite[Lemma 3.8]{fu}. We note that 
we can apply \cite[Lemma 3.8]{fu} by the flattening theorem. 
We note that 
$$
K_X = g^*(K_Z+L_{X/Z}^{ss})+
\sum_{D_i} s_{D_i} g^*D_i
+ B, 
$$ 
where $D_i$ is an irreducible component of $D$ for every $i$ 
(see Section \ref{3}). 
\end{itemize}

By the canonical bundle formula, 
we have 
$$
g_{*} K_{X/Z}^{\otimes m}(mB_{-})\simeq {\mathcal O}_{Z}
(\sum _{i}ms_{D_{i}}D_{i}),
$$ 
where $m$ is a positive integer such that 
$ms_{D_{i}}$ are integers for every $i$. 
We note that we can assume that 
$mL_{X/Z}^{ss}$ is trivial by Theorem \ref{dou} since $\xVar(g)=0$. 
By restricting the canonical bundle formula to 
$X_{y}\to Z_{y}$, where $y$ is a sufficiently 
general point of $Y$, 
we obtain an irreducible component $D_{0}$ of $D$ such that 
$h_1\nu(D_0)=Y$ and $s_{D_{0}}\ne 0$ since $\kappa 
(X_{y})=\dim Z_y\geq 1$. 

Let $\overline{D}_{0}$ be the image of $D_{0}$ on $Y\times A$. 
Then $\kappa (Y\times A, \overline{D}_{0})>0$ by Corollary \ref{new} 
below. 
On the other hand, 
every irreducible component of $\nu^{*}\overline{D}_{0}-D_{0}$ is 
$\nu$-exceptional and 
$$H^{0}(Z,{\mathcal O}_{Z}
(ms_{D_{0}}(D_{0}-\nu^{*}{\overline{D}_{0}}))
\otimes K_{Z/Y}^{\otimes mk})\ne 0$$ for a sufficiently large integer
$k$. 
We note that $K_{Y\times A}=h_{1}^{*}K_{Y}$ since 
$A$ is an Abelian variety. 
Combining the above, 
we obtain 
$$H^{0}(Z,g_{*}K_{X/Y}^{\otimes mk}(kmB_{-})\otimes{\mathcal O}_{Z}
(-ms_{D_{0}}\nu ^{*}{\overline{D}_{0}}))\ne 0.$$ 
Therefore, 
$$\kappa (X, K_{{X} /Y})
\geq \kappa (Z, \nu ^{*}{\overline{D}_{0}})=
\kappa (Y\times A, \overline{D}_{0})>0.$$
We note that $B_-$ is effective and exceptional over $\widetilde X$. 
Thus, we finish the proof. 
\endproof

\section{Some remarks on Abelian varieties}\label{5}

The main purpose of this section is 
to prove Corollary \ref{new}, which was already used in 
the proof of the main theorem. The results below are variants of 
the theorem of cube. 

\begin{say}
Let $Y$ be a variety ($Y$ is not necessarily complete) 
and $A$ an Abelian variety. 
We define $Z:=Y\times A$. 
Let $\mu:A\times A\to A$ be the multiplication. 
Then $A$ acts on 
$A$ naturally by the group law of $A$. 
This action induces a natural action on $Z$. 
We write it by $m:Z\times A\to Z$, that is, 
$$
m:((y,a),b)\to (y,a+b),
$$ 
where $(y,a)\in Y\times A=Z$ and $b\in A$.
Let $p_{1i}:Z\times A\times A\to Z\times A$ 
be the projection onto the $(1,i)$-th factor for $i=2,3$ 
and $p_{23}:Z\times A\times A\to A\times A$ the 
projection onto the $(2,3)$-factor. 
Let $p:Z\times A\times A\to Z$ be the first projection 
and $p_i:Z\times A\times A\to A$ the 
projection onto the $i$-th factor for $i=2,3$. 
We define the projection $\rho:Z=Y\times A \to A$. 
We fix a section $s:A\to Z$ such that $s(A)=\{y_0\}\times 
A$ for a point $y_0\in Y$. 
We define the morphisms as follows; 
\begin{eqnarray*}
\pi_i&:=&p_i\circ(s\times id_A \times id_A)
\ \  \text{for}\ \  i=2,3,\\
\pi_{23}&:=&p_{23}\circ(s\times id_A \times id_A), \\
\pi&:=&\rho\times id_A\times id_A. 
\end{eqnarray*}
Let $L$ be a line bundle on $Z$. 
We define a line bundle $\mathcal L$ 
on $Z\times A\times A$ as follows; 
\begin{eqnarray*}
\mathcal L=&p^*L\otimes (id_Z \times \mu)^{*}m^{*}L
\otimes (p_{12}^{*}m^{*}L)^{-1}\otimes 
(p_{13}^{*}m^{*}L)^{-1}\\ 
&\otimes \pi^{*}((\pi_{23}^{*}\mu^{*}s^{*}L)^{-1}
\otimes \pi_2^{*}s^{*}L\otimes \pi_3^{*}s^{*}L). 
\end{eqnarray*}
\end{say}

\begin{thm}\label{cubic}
Under the above notation, 
we have that 
$$
\mathcal L\simeq \mathcal O_{Z\times A\times A}. 
$$
\end{thm}
\proof
It is not difficult to check that 
the restrictions $\mathcal L$ to 
each of $Z\times \{0\}\times A$ and 
$Z\times A\times \{0\}$ are 
trivial by the definition of $\mathcal L$, where $0$ is the origin of $A$. 
We can also check that the 
restriction onto $s(A)\times A\times A$ is trivial 
(cf.~\cite[p.58 Corollary 2]{Mu}). 
In particular, $\mathcal L|_{\{z_0\}\times A\times A}$ is 
trivial for any point $z_0\in s(A)\subset Z$. 
Therefore, by the theorem of cube \cite[p.55 Theorem]{Mu}, 
we obtain that $\mathcal L$ is trivial.   
\endproof

We write $T_a:=m|_{Z\times \{a\}}:Z\simeq Z\times \{a\}\to Z$, 
that is,  
$$
T_a:(y,b)\to (y, b+a),  
$$ 
for $(y,b)\in Y\times A=Z$. 

\begin{cor}\label{cube2}
By restricting $\mathcal L$ to $Z\times \{a\}\times\{b\}$, 
we obtain; 
$$
L\otimes T_{a+b}^{*}L\simeq T_{a}^{*}L\otimes 
T_{b}^{*}L, 
$$ 
where $a, b\in A$. 
\end{cor}

The following is a supplement and a generalization of \cite[Lemma 7.11]{fu}. 

\begin{cor}\label{new}
Let $D$ be a Cartier divisor on $Z$. 
Then $2D\sim T_a^{*}D+T_{-a}^{*}D$ for 
$a\in A$. 
In particular, if $Y$ is complete and 
$D$ is effective and not vertical with 
respect to $Y\times A\to Y$, then 
$\kappa (Z, D)>0$.   
\end{cor}
\proof 
We put $L=\mathcal O_Z(D)$ and $y=-a$. 
Apply Corollary \ref{cube2}. 
We note that $\Supp D\ne\Supp T_a^{*}D$ if 
we choose $a\in A$ suitably. 
\endproof

\ifx\undefined\bysame
\newcommand{\bysame}{\leavevmode\hbox to3em{\hrulefill}\,}
\fi

\end{document}